\documentclass[12pt,english]{article}
\usepackage[T1]{fontenc}
\usepackage[latin9]{inputenc}
\usepackage{setspace}
\usepackage{babel}
\usepackage{amsfonts,mathrsfs}
\usepackage{amssymb,color}
\usepackage{booktabs}

\usepackage{graphicx}
\usepackage{epsfig}

\usepackage{fancyhdr} 
\newtheorem{theorem}{Theorem}
\newtheorem{corollary}{Corollary}

\def\qed{\quad{$\square$}}
{\bf}{}

\title{A game-theoretic approach to dynamic demand response management\thanks{This work  has been supported by the 2012 ``Research Fellow'' Program of the Dipartimento di Matematica, Universit\`a di Trento and by the PRIN 20103S5RN3 ``Robust decision making in markets and organization".}
} 

\author{D. Bauso\thanks{Dario Bauso is with Dipartimento di Ingegneria Chimica, Gestionale, Informatica e Meccanica,
Universit\`a di Palermo, Italy, and Dipartimento di Matematica, Universit\`a di Trento (as research fellow), Via Sommarive 14, 38050 Povo-Trento, Italy.   D. Bauso is currently academic visitor at the Department of Engineering Science, University of Oxford, United Kingdom.
email: \textsl{dario.bauso@unipa.it}} 
}

\begin{document}

\maketitle
\thispagestyle{plain}
\pagestyle{plain}

\begin{abstract}                          

Within the realm of dynamic of  \emph{smart buildings} and \emph{smart cities}, dynamic response management is playing an ever-increasing role thus attracting the attention of scientists from different disciplines. Dynamic demand response management involves a set of operations aiming at decentralizing the control of loads in  large and complex power networks. Each single appliance if fully responsive and readjusts its energy demand to the overall network load. A main issue is related to mains frequency oscillations resulting from an unbalance between supply and demand. In a nutshell, this paper contributes to the topic by equipping each signal consumer with strategic insight. In particular, we highlight three main contributions and a few other minor contributions. First, we design a mean-field game for the TCLs application, study the mean-field equilibrium for the deterministic mean-field game and investigate on asymptotic stability for the microscopic dynamics. Second, we extend the analysis and design to imperfect models which involve both stochastic or deterministic disturbances. This leads to robust mean-field equilibrium strategies guaranteeing stochastic and worst-case stability, respectively. Minor contributions involve the use of stochastic control strategies rather than deterministic, and some numerical studies illustrating the efficacy of the proposed strategies. 

\end{abstract}


\section{Introduction}

Demand response management involves a set of operations aiming at decentralizing load control in  power networks \cite{AS07,ENT07,GC88,US06}. In particular, it calls for the alteration of the timing, of the level of
instantaneous demand, or of the total electricity by end-use customers from their normal consumption patterns in response to changes in the price of electricity over time. This is possible also through an opportune design of  incentive payments  to induce lower electricity use at off-peak times.

\begin{figure} [h]
\centering
\def\svgwidth{.55\columnwidth}
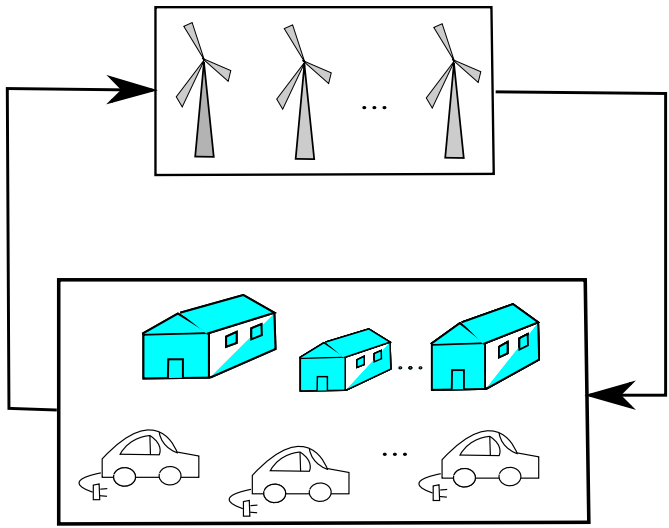
\caption{Demand response involves populations of electrical loads (lower block) and energy generators (upper block) intertwined in a feedback-loop scheme.}  \label{fig:demandresponse}
\end{figure}

A communication protocol aggregates relevant information on the past, current and forecasted demand and transmits it to each fully responsive load controller or decision mechanism, which will adopt opportune actions such as increasing or decreasing the proper load or energy  demand. The novelty of this paper is in that \emph{fully responsive load control} together with the many cooperative and competitive aspects involved in the process, are now reviewed as a game with a large number of indistinguishable players, these being the single loads. For illustrative purposes, in this paper, fully responsive load control is reviewed in the context of thermostatically controlled loads (TCLs),  in  smart buildings or plug-in electric vehicles \cite{AK12,MCH13,MKC13,PCGL14,bcSGEFA14}, see Fig. \ref{fig:demandresponse}.

A first idea of this work, which is common also to \cite{AK12,BagBau13}, is to adopt stochastic response strategies rather than deterministic. This means that each TCL selects a probability with which to activate one of the two functioning modes, \emph{on} and \emph{off}. Thus a probability value of $1/2$ means that the TCL is $50 \%$ \emph{on} and $50 \%$ \emph{off}. It has been shown in \cite{AK12,BagBau13} that stochastic response strategies outperform deterministic ones, especially in terms of  attenuating  the \emph{mains frequency} oscillations. These are due to the unbalance  between energy demand and supply (see e.g. \cite{RDM12}). The mains frequency usually needs to be stabilized around a nominal value (50 Hz in Europe). If electrical demand exceeds generation then frequency will decline, and vice versa. 

\begin{figure} [htb]
\centering
\includegraphics[width=\columnwidth]{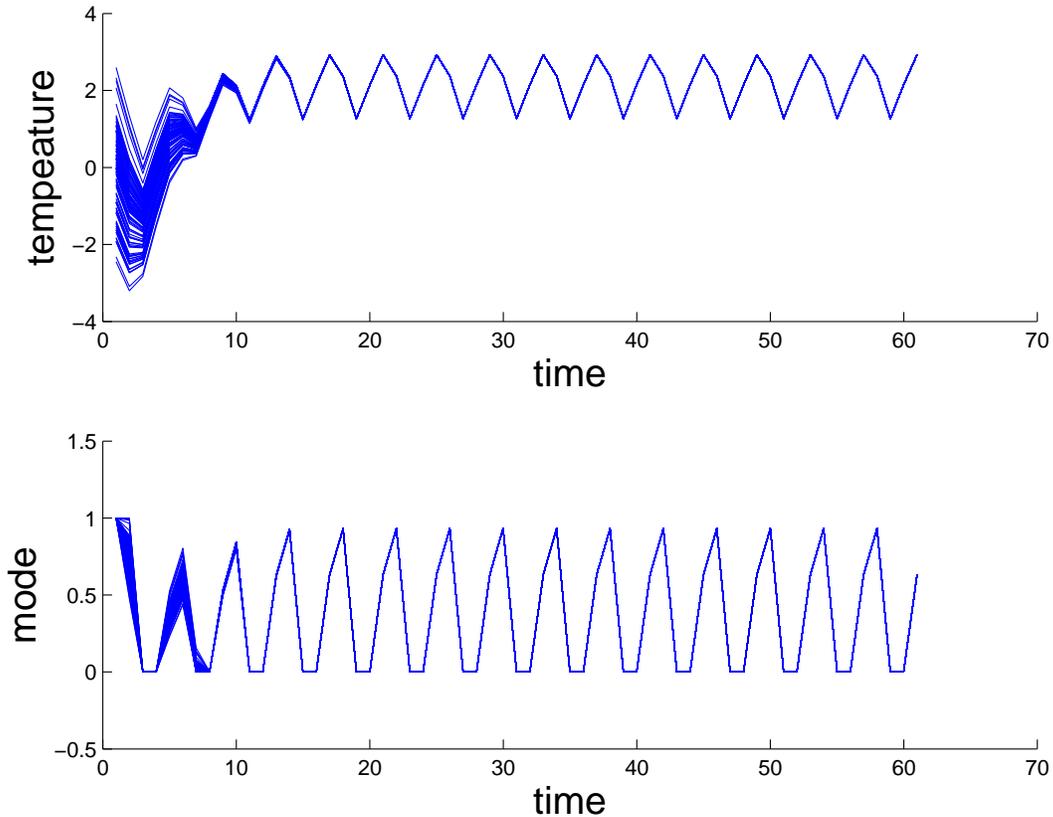}
\caption{Example of oscillations: qualitative time plot  of the state of each TCL, namely temperature  (top row) and mode of functioning (bottom row).}
\label{fig:oscillations}
\end{figure}

A qualitative plot of such an oscillatory phenomenon is displayed in Fig. \ref{fig:oscillations}. The two rows depict the time plot  of the state of each TCL, namely the temperature in the top row and the mode of functioning in the bottom row. Here each TCL increases or decreases its proper load in response to the current network load, and as clear visually, this  induces oscillations in the mains frequency due to an undesired synchronized reaction of the whole population of TCLs. 
 
This preamble introduces the main aim of this paper, which studies constructive design methods of distributed demand response management strategies in order to reduce the mains frequency oscillations and stabilize both the temperature and the functioning mode of the TCLs.  

The model used in this paper is as follows. Each single TCL is a player and is characterized by two state variables, the temperature and the functioning mode. The state dynamics of a TCL  --- henceforth referred to as \emph{microscopic} dynamics to distinguish it from the dynamics of the aggregate temperature and functioning mode of the whole population, the latter called \emph{macroscopic} dynamics --- describes the time evolution of its temperature and mode in the form of a linear ordinary differential equation in the deterministic case, and of a stochastic differential equation in the stochastic case. In addition to the state dynamics, each TCL is programmed with a given finite-horizon cost functional that accounts for i)  energy consumption, ii) deviation of mains frequency from the nominal one, and iii) deviation of the TCL's temperature from a reference value.  More formally, the mains frequency involved in the specifics ii) is used in a cross-coupling mean-field term that incentivizes the TCL to switch to $off$ if the mains frequency is below the nominal value and to switch to $on$ if the mains frequency is above the nominal value. In other words, the cross-coupling mean-field term models all kinds of incentive payments, benefits, or smart pricing policies aiming at shifting demand from high-peak to off-peak periods.

\subsection{Highlights of contributions}
 This paper provides three main results. First, in the spirit of prescriptive game theory and mechanism design \cite{BagBau} we design a mean-field game for the TCLs application, study the  mean-field equilibrium for the deterministic mean-field game and investigate on asymptotic stability for the microscopic dynamics. Asymptotic stability means  that both the temperature and the mode functioning of each TCL converges to the reference value. A second result relates to the stochastic case, characterized by a stochastic disturbance in the form of a Brownian motion in the microscopic dynamics. After establishing a mean-field equilibrium, we provide some results on stochastic stability. In particular, we focus on two distinct scenarios. In one case, we assume that the stochastic disturbance expires in a neighborhood of the origin. This reflects in having the Brownian motion coefficients linear in the state. The resulting dynamics is well-known in the literature as geometric Brownian motion. As for any geometric Brownian motion, we can study conditions for it to be stochastically  stable almost surely. This means that the state trajectories converge to zero with probability one.  In a second case, the stochastic disturbance is independent on the state and the Brownian motion coefficients are constant. This leads to a  dynamics which resembles  the Langevin equation. Following well-known results on the Langevin equation,  the dynamics is proven to be stochastically stable in the second-moment. An expository work on stochastic analysis and stability is \cite{LF96}. A third result deals with robustness for the microscopic dynamics. The dynamics is now influenced by an additional adversarial disturbance, with bounded resource or energy. Even for this case, we study the mean-field equilibrium and  investigate on conditions that guarantee worst-case stability.

\subsection{Literature overview}
We introduce next two streams of literature. One is related to dynamic response management, while the second one is about the theory of differential games with a large number of indistinguishable players, also known as mean-field games. 

\subsubsection{Related literature on demand response}  

Examples of papers developing the idea of dynamic demand management  are \cite{CH11,CPTD12,MCH13,MKC13}. In particular, \cite{CH11} provides an overview on  the redistribution of the load away from peak hours and the design of decentralized strategies to produce a predefined load trajectory. This idea is further developed in  \cite{CPTD12}. 
  To understand the role of game theory in respect to this specific context the reader is referred to \cite{MCH13}. There, the authors present a large population game where the agents are plug-in electric vehicles and the Nash-equilibrium strategies (see \cite{B99}) correspond to distributed charging policies that redistribute the load away from peaks. The resulting strategies are known with the name of  \emph{valley-filling strategies}. In this paper we adopt the same perspective in that we show that network frequency stabilization can be achieved by giving incentives to the agents to adjust their strategies in order to converge to a mean field equilibrium. To do this, in the spirit of \emph{prescriptive game theory} \cite{BagBau}, a central planner or game designer has to design the individual objective function so to penalize those agents that are in $on$ state in peak hours, as well as those who are in $off$ state in off-peak hours.  Valley-filling and coordination strategies have been shown particularly efficient in thermostatically controlled loads such as refrigerators, air conditioners and electric water heaters \cite{MKC13}. 


The results obtained in this paper are in accordance with the recent
results in \cite{AK12}, according to which, stochastic control laws are in general more appropriate than deterministic ones when it comes to desynchronize the appliances functioning. 


\subsubsection{Related literature on mean-field games}
 A second stream of literature related to the problem at hand is on mean-field games.  
Mean-field games were formulated by  Lasry and Lions in \cite{LL07} and  independently by 
M.Y. Huang, P. E. Caines and R. Malham\'e  in \cite{HCM06,HCM07}.  The mean-field theory of dynamical games is a modeling framework at the interface of differential game theory, mathematical physics, and $H_\infty$-optimal control that tries to capture the mutual influence between a crowd and its individuals. 

From a mathematical point of view the mean-field approach leads to a system of two PDEs. The first PDE is the Hamilton-Jacobi-Bellman (HJB) equation. 
The second PDE is the Fokker-Planck-Kolmogorov (FPK) equation which describes the density of the players. 
Explicit solutions in terms of mean-field equilibria are available for  linear-quadratic mean-field games \cite{B12}, and have been recently extended to more general cases in \cite{Gomes:14}. In addition to explicit solutions, a variety of solution schemes have been recently proposed based on discretization and/or numerical approximations, see e.g. \cite{ACC12,AC10,PB11}. The idea of extending the state space, which originates in optimal control \cite{SA12,SA13}, has been also used to approximate mean-field equilibria in \cite{BMA-ECC14}. 
More recently, robustness and risk-sensitivity have been brought into the picture of mean-field games \cite{BTB12,TZB14}, where the first PDE is now the Hamilton-Jacobi-Isaacs (HJI) equation. For a survey on mean-field games and applications we refers the reader to \cite{GLL10}. A first attempt to apply mean-field games to demand management is in \cite{BagBau13}. 

The paper is organized as follows. In Section \ref{sec:model} we state the problem and introduce the model.
In Section \ref{sec:prel} we review some preliminary results. In Section \ref{mainresults} we state and discuss the main results. In Section \ref{discussion} we provide some discussion. In Section \ref{sec:num} we carry out some numerical studies. Finally, in Section \ref{conclusion} we provide some conclusions.


\subsection{Notation}

The symbol $\mathbb E$ indicates the expectation operator. We use $\partial_x$ and $\partial^2_{xx}$ to denote the first and second partial derivatives with  respect to $x$, respectively. Given a vector $x \in \mathbb R^n$ and a matrix $a\in \mathbb R^{n \times n}$ we denote by $\|x\|^2_{a}$ the weighted two-norm $x^T a x$. The symbol $a_{i \bullet}$ means the $i$th row of a given matrix $a$. We denote by $Diag(x)$ the diagonal matrix in $\mathbb R^{n \times n}$ whose entries in the main diagonal are the components of $x$.  We denote by $dist(X,X^*)$ the distance between two points $X$ and $X^*$ in $\mathbb R^n$. We denote by $\Pi_{\mathcal M}(X)$ the projection of $X$ onto set $\mathcal M$. The symbol ``:'' denotes the Frobenius product.

\section{Population of TCLs through mean-field games}\label{sec:model}
In this section, in the spirit of prescriptive game theory and mechanism design \cite{BagBau}, we design a mean-field game for the TCLs application, with the aim of incentivizing cooperation among the TCLs through an opportune design of distributed cost functionals, one per each TCL. 

Consider a population of hybrid controlled thermostat loads (TCLs) and a time horizon window $[0,T]$. Each TCL is characterized by a continuous state, namely the temperature $x(t)$, and a binary state $\pi_{on}(t) \in \{0,1\}$, representing the condition $on$ or $off$ at time $t\in [0,T]$.  
 When the TCL is set to $on$ the temperature decreases exponentially up to a fixed lower temperature $x_{on}$ whereas in the $off$  position the temperature increases exponentially up to a higher temperature $x_{off}$. Then, the temperature of each appliance evolves according to the following differential equations: 
\begin{equation}\label{dyn0} 
\dot x(t)= \left\{ \begin{array}{ll}  -\alpha (x(t)-x_{on}) & \mbox{if $\pi_{on}(t)=1$}\\
 -\beta (x(t)-x_{off}) & \mbox{if $\pi_{on}(t)=0$}  \end{array}, \ t\in [0,T), \right. \end{equation}
with initial state $x(0)=x$ and where the rates $\alpha,\beta$ are given positive scalars. 
\bigskip
 
\begin{figure} [h]
\centering
\def\svgwidth{.45\columnwidth}
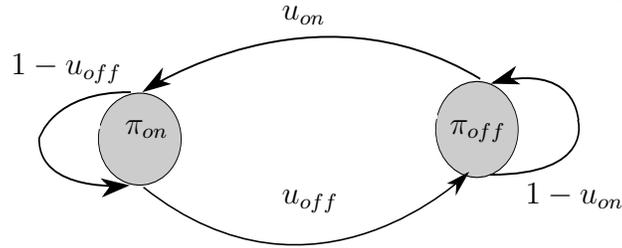
\caption{Automata describing transition rates from $on$ to $off$ and viceversa.}  \label{fig:transition}
\end{figure}

In accordance with \cite{AK12,BagBau13} we set the problem in a stochastic framework where each TCL is in one of the two states $on$ or $off$ 
with  given probabilities  $\pi_{on} \in [0,1] $ and $\pi_{off} \in [0,1]$. The control variable is the transitioning rate 
$u_{on}$ from $off$ to $on$ and the transitioning rate  $u_{off}$ from $on$ to $off$. This is illustrated in the automata in 
Fig. \ref{fig:transition}. 

The corresponding dynamics is then given by
\begin{equation}\label{dynoo} 
\left\{ \begin{array}{lll}  \dot \pi_{on}(t)  =  u_{on} (t) - u_{off}(t) , \ t\in [0,T),\\
\dot \pi_{off}(t)   =  u_{off} (t) - u_{on}(t) , \ t\in [0,T),\\
0\leq  \pi_{on}(t),\pi_{off}(t) \leq 1,\ t\in [0,T).
\end{array} \right. \end{equation}

As $\dot \pi_{on}(t) + \dot \pi_{off}(t)=0$, we can simply consider only one of the above dynamics. Then, let us denote $y(t)=\pi_{on}(t)$ and introduce a stochastic disturbance in the form of a Brownian motion, denote it $\mathcal B(t)$, and a deterministic disturbance $w(t)=[w_1(t)\quad w_2]^T$. For any $x,y$ in the $$\mbox{``set of feasible states''} \quad \mathcal S:=]x_{on},x_{off}[ \times ]0,1[,$$ the resulting dynamics in a very general form is given by 
\begin{equation}\label{dyno} 
\left\{ \begin{array}{lll} d x(t) & = & \Big(y(t) \Big[ -\alpha (x(t)-x_{on})\Big] + (1-y(t))\Big[ -\beta (x(t)-x_{off})\Big]\\ &&+ d_{11} w_1(t) +d_{12} w_2(t)\Big) dt+  \sigma_{11}(x) d\mathcal B(t),\\
& =: &\Big( f(x(t),y(t))+ d_{11} w_1(t) +d_{12} w_2(t)\Big)dt+ \sigma_{11}(x) d\mathcal B(t) , \ t\in [0,T),\\
x(0)& = & x,\\
d y(t) & = & \Big(u_{on}(t) - u_{off}(t)+ d_{21} w_1(t) +d_{22} w_2(t)\Big)dt + \sigma_{22}(y) d\mathcal B(t)  \\
& =: & \Big(g(u(t))+ d_{21} w_1(t) +d_{22} w_2(t)\Big)dt + \sigma_{2}(y) d\mathcal B(t),\ t\in [0,T),\\
y(0) & = & y, \end{array} \right. \end{equation}
where $\sigma_{ij}$ and $d_{ij}$, $i,j=1,2$ are positive scalar coefficients. 

For a mean-field game formulation, consider a  probability density function $m:[x_{on},x_{off}]\times [0,1] \times[t,T] \to[0,+\infty[$, $(x,y,t) \mapsto m(x,y,t)$,  which satisfies $\int_{x_{on}}^{x_{off}} \int_{[0,1]} m(x,y,t)dxdy=1$ for every $t$.  Let us also define 
as $m_{on}(t) := \int_{x_{on}}^{x_{off}} \int_{[0,1]} y m(x,y,t)dxdy$. Likewise we denote by $m_{off}(t)=1-m_{on}(t)$.

At every time $t$ the network or mains frequency depends linearly on the discrepancy  between the percentage of TCLs in $on$ position and a nominal value. We call such a discrepancy as error and denote it by $e(t)=m_{on}(t) -\overline m_{on}$, where  $\overline m_{on}$ is the nominal value (the higher the percentage of TCLs in $on$ position with respect to the nominal value, the lower the network frequency).   

We then consider the running cost below, which depends on the distribution $m(x,y,t)$ through the error $e(t)$:
\begin{equation}\label{g}
\begin{array}{r}  c(x(t),y(t),u(t),m(x,y,t))  =     \frac{1}{2} \Big(q x(t)^2 + r_{on} u_{on}(t)^2+ r_{off} u_{off}(t)^2\Big) \\+ y(t) (S e(t)  + W),
\end{array}
\end{equation}
where $q$, $r_{on}$, $r_{off}$, and $S$ are opportune positive scalars.

Note that cost (\ref{g}) includes four terms. The term $\frac{1}{2}q x(t)^2$ penalizes  the deviation of the TCLs' temperature from the nominal value, which we set to zero. The terms  $\frac{1}{2} r_{on} u_{on}(t)^2$  introduces a cost for fast switching; namely this cost is zero when either $u_{on}(t)=0$ (no switching) and is maximal when $u_{on}(t)=1$ (probability 1 of switching). Similar comment applies $\frac{1}{2} r_{off} u_{off}(t)^2$. The term $y(t) S e(t)$ accounts for the network stabilization in that it penalizes those appliances that are $on$ whenever $e(t)>0$, the latter condition meaning that demand exceeds supply. The same term turns into a revenue if an appliance is $on$ whenever $e(t)< 0$, i.e., whenever supply exceeds demand. Finally, the penalty term $y(t) W$  accounts for the minimization of power, namely, whenever the TCL is $on$ the power consumption is $W$.

Also consider a terminal cost $\Psi:\mathbb{R}\to[0,+\infty[$, $x\mapsto\Psi(x)$ to be yet designed. 
\bigskip

\noindent 
 \textbf{Problem statement}. Given a finite horizon $T>0$ and an initial distribution $m_0:[x_{on},x_{off}] \times [0,1] \to [0,+\infty[$, minimize over $\mathcal U$ and maximize over $\mathcal W$ , subject to the controlled system (\ref{dyno}), the cost functional
\[
J(x,y,t,u(\cdot))=\mathbb E \int_0^T (c(x(t),y(t),u(t),w(t),m(x,y,t))- \frac{1}{2} \gamma^2 \|w(t)\|^2 )d t + \Psi(X(T)),
\]
\noindent
where $\gamma$ is a positive scalar,  $\mathcal{U}$ and $\mathcal{W}$ are the sets of all measurable state feedback closed-loop policies  $u(\cdot):[0,+\infty[\to \mathbb R$ respectively, and $w(\cdot):[0,+\infty[\to \mathbb R$
and $m(\cdot)$ is the time-dependent function describing the evolution of the mean of the  distribution of the TCLs' states.

\section{Preliminary results}\label{sec:prel}
This section reviews first- and second-order mean-field games in preparation to apply the game to the problem at hand. In the first case, the microscopic dynamics is deterministic and the resulting mean-field game involves only the first derivatives of the value function and of the density function. In the second case, the microscopic dynamics is a stochastic differential equation driven by a Brownian motion, which leads to the involvement of second derivatives of the value function and density function.  In addition to this, this section  specializes the model to the case under study, involving a population of TCLs and introduced in the previous section.  
\subsection{First- and second-order mean-field games}
This section streamlines some preliminary results on mean-field games. To this purpose, consider a generic cost and dynamics 
\begin{equation}
\begin{array}{cll}J(X,0,U(.))=\inf_{U(.)} \int_{t=0}^T c(X(t),m,U(.)) d t + \Psi(X(T)), \\ \\
\dot X(t) = F(X(t),U(.))  \mbox{ in } \mathbb R^n,
\end{array}\end{equation}
where $c(.)$ is the running cost,  $\Psi(X)\ \forall\ X\in \mbox{in } \mathbb R^n$ is  the terminal penalty, and
where $U(.)$ is any state-feedback closed loop control policy.
Let $v(X,t)$ be the value function, i.e., the optimal value of $J(X,t,U(\cdot))$. Then from \cite{LL07} it is well-known that the problem results in the following mean-field game system

\begin{equation}
\label{eq:meanfieldprel}
\left\{
\begin{array}{ll}
\displaystyle
-\partial_t v(X,t)-F(X,U^*(X)) \partial_X v(X,t) - c(X,m,U^*(X))=0 & \mbox{(a)} \\
 \mbox{in } \mathbb R^n \times]0,T],\\
\displaystyle
v(X,T)=\Psi(X)\ \forall\ X\in \mbox{in } \mathbb R^n,\\\\
\displaystyle 
U^*(X,t)=\mbox{argmax}_{U\in \mathbb R}\{-F(X,U) \partial_X v(X,t) - c(X,m,U)\},& \mbox{(b)}
\end{array}
\right.
\end{equation}

\begin{equation}
\label{eq:meanfieldprelm}
\left\{
\begin{array}{ll}
\displaystyle
\partial_t m(X,t)+ div (F(X,U^*(X)) m(X,t))=0\ 
 \mbox{in } \mathbb R^n \times]0,T], \\
 m(X,0) = m_0(X), \forall\ X\in \mbox{in } \mathbb R^n.
\end{array}
\right.
\end{equation}
The partial differential equation (PDE) \ref{eq:meanfieldprel} (a) is the Hamilton-Jacobi-Bellman equation which returns the value function $v(X,t)$ once we fix the distribution $m(X,t)$; This PDE has to be solved backwards with boundary conditions at final time $T$, represented by the last line in \ref{eq:meanfieldprel} (a). In \ref{eq:meanfieldprel} (b) we have the optimal closed-loop control $U^*(X,t)$ as maximizer of the Hamiltonian function in the rhs. The PDE \ref{eq:meanfieldprelm} represents the transport equation of the measure $m$ immersed in a vector field $F(X,U^*(X))$; It returns the distribution $m(X,t)$ once fixed the  the optimal closed-loop control $U^*(X,t)$ and consequently the vector field $F(X,U^*(X))$. Such a PDE has to be solved forwards with boundary condition at the initial time (see the last line) of (\ref{eq:meanfieldprelm}). 

In a second order mean-field game, the dynamics is a stochastic differential equation driven by a Brownian motion, and the cost function is considered through its expected value, namely,
\begin{equation}
\begin{array}{cll}J(X,0,U(.))= \inf_{U(.)} \mathbb E \int_{t=0}^T c(X(t),m,U(X(t))) d t + \Psi(X(T)) \\ \\
d X(t) = F(X(t),U(.)) dt + \sigma(X) d\mathcal B(t) \mbox{ in } \mathbb R^n,
\end{array}\end{equation}
where $\mathcal B(t) \in \mathbb R^n$ is the Brownian motion and  $\sigma(X) \in \mathbb R^{n \times n}$ is the coefficient matrix.

From \cite{LL07} the second-order mean-field game system is then given by
\begin{equation}
\label{eq:meanfieldprel2}
\left\{
\begin{array}{ll}
\displaystyle
-\partial_t v(X,t)-F(X,U^*(X)) \partial_X v(X,t) - c(X,m,U^*(X))\\
 - \frac{1}{2}  \sigma(X)\sigma(X)^T : \partial_{XX}  v(X,t) =0 
 \mbox{in } \mathbb R^n \times]0,T],& \mbox{(a)}\\
\displaystyle
v(X,T)=\Psi(X)\ \forall\ X\in \mbox{in } \mathbb R^n,\\\\
\displaystyle 
U^*(X,t)=\mbox{argmax}_{U\in \mathbb R}\{-F(X,U) \partial_X v(X,t) - c(X,m,U)\},& \mbox{(b)}
\end{array}
\right.
\end{equation}
\begin{equation}
\label{eq:meanfieldprel2m}
\left\{
\begin{array}{ll}
\displaystyle
\partial_t m(X,t)+ div (F(X,U^*(X)) m(X,t)) \\- \frac{1}{2} \sum_{i=1}^n \sum_{j=1}^n \partial^2_{X_iX_j}( \tilde \sigma_{ij} m(X,t)) =0\ 
 \mbox{in } \mathbb R^n \times]0,T],\\
 m(X,0) = m_0(X), \forall\ X\in \mbox{in } \mathbb R^n,
\end{array}
\right.
\end{equation}
where the symbol ``:'' denotes the Frobenius product and  $\tilde \sigma_{ij} = \sum_{k=1}^n \sigma_{ik}(X) \sigma_{jk}(X)$. 

In a second-order mean-field game the Hamilton-Jacobi-Bellman equation, as in \ref{eq:meanfieldprel2} (a),  involves the second-order derivatives of the value function in the additional term represented by the Frobenius product;  Likewise, also the transport equation as in (\ref{eq:meanfieldprel2m}) involves the second-order derivatives of the density function. The rest of the system is similar to the first-order case. Let us now specialize the above model to the TCLs application introduced in the previous section.  

\subsection{Mean-field game for the TCL application}

Specializing to our TCLs application, let $v(x,y,m,t)$ be the value function, i.e., the optimal value of $J(x,y,t,u(\cdot))$. 
 Let us denote by
 $$\nonumber\label{coeff} 
k(x(t))  =  x(t) (\beta - \alpha) + ( \alpha x_{on} - \beta x_{off}).
$$
Then,  the problem at hand can be rewritten as
 in terms of the state, control and disturbance vectors 
$$X(t)=\left[\begin{array}{ll}
x(t) \\
y(t)\end{array}\right],\quad u(t)=\left[\begin{array}{lll} u_{on}(t)\\ u_{off}(t)\end{array}\right], \quad w(t)=\left[\begin{array}{lll} w_1(t)\\ w_2(t)\end{array}\right] 
$$
and yields the   linear quadratic problem:
\begin{equation}\label{pq1}
\begin{array}{cll}
\displaystyle
\inf_{\{u_t\}_{t}}  \mathbb E \int_0^T \left[\frac{1}{2} \Big(\|X(t)\|^2_Q  +  \|u(t)\|^2_R  -\gamma^2 \|w(t)\|^2\Big)+ L^T X(t) \right] dt\\ \\
d X(t) = (A X(t) +  B u(t) + C +Dw(t)) dt + \Sigma d\mathcal B(t),   \mbox{ in } \mathcal S
\end{array}
\end{equation}
where 
\begin{equation}\label{fd1}\nonumber
\begin{array}{ccc}
Q = \left[\begin{array}{ccc}
q & 0\\
0 &  0 \end{array}\right], & R =r=\left[\begin{array}{cc} r_{on} & 0 \\ 0 & r_{off} \end{array}\right], &
 L(e)=\left[\begin{array}{cc} 0 \\ S e(t) + W \end{array}\right],\\ \\
 A(x) = \left[\begin{array}{ccc}
-\beta & k(x(t)) \\
0 & 0 \end{array}\right],
& 
B=\left[\begin{array}{cc}
0 & 0\\
1 & -1 \end{array}\right], & C = \left[\begin{array}{c}
\beta x_{off} \\
0
\end{array}\right],
\end{array}
\end{equation}
and 
\begin{equation}\label{fd2}\nonumber
\begin{array}{cccc}
D = \left[\begin{array}{ccc}
d_{11} & d_{12}\\
d_{21} &  d_{22} \end{array}\right],  \quad 
  \Sigma = \left[\begin{array}{ccc}
\sigma_{11}(x) & 0\\
0 &  \sigma_{22}(y) \end{array}\right] 
\end{array}.
\end{equation}

The resulting mean-field game is given by
\begin{equation}\label{mainHJI}\nonumber
\left\{
\begin{array}{lll}
\partial_t \mathcal V_t(X) + \inf_u \sup_w\Big\{  \partial_X \mathcal V_t(X)^T (AX + Bu +C+D w) + \frac{1}{2} \Big(\|X\|^2_Q \\
 \qquad \qquad +  \|u\|^2_R  - \gamma^2 \|w\|  \Big)+ L^T X \Big\}+ \frac{1}{2} ( \sigma_{11}(x)^2  \partial_{xx} v(X,t),& \mbox{(a)} \\  \qquad \qquad  \qquad 
 +\sigma_{22}(y)^2  \partial_{yy} v(X,t))=0,   \mbox{ in } \mathcal S \times [0,T[,\\ 
v(X,T)=g(x),   \mbox{ in } \mathcal S
\\\\
\displaystyle 
u^*(x,t)=\mbox{argmin}_{u\in \mathbb R} \Big\{  \partial_X \mathcal V_t(X)^T (AX + Bu +C+D w) + \frac{1}{2} \|u(t)\|_R^2  \Big\},
& \mbox{(b)}\\
\displaystyle 
w^*(x,t)=\mbox{argmax}_{u\in \mathbb R} \Big\{  \partial_X \mathcal V_t(X)^T (AX + Bu +C+D w) - \frac{1}{2}  \gamma^2 \|w(t)\|^2 \Big\}
\end{array}
\right.
\end{equation}
%
and 
\begin{equation}
\label{mainFPK}
\left\{
\begin{array}{ll}
\displaystyle
\partial_t m(x,y,t)+ div [(AX + Bu +C+D w) \ m(x,y,t)] \\
- \frac{1}{2} \sum_{i=1}^2 \sum_{j=1}^2 \partial^2_{X_iX_j}( \tilde \sigma_{ij} m(X,t)) =0\ 
 \mbox{ in } \mathcal S \times]0,T[,\\
\displaystyle
m(x_{on},y,t)=m(x_{off},y,t)=0\ \forall\ y\in[0,1],\ t\in[0,T],\\
\displaystyle
m(x,y,0)=m_0(x,y)\ \ \forall\ x\in [x_{on},x_{off}],\ y \in [0,1]\\
\displaystyle
\int_{x_{on}}^{x_{off}}m(x,t)dx=1\ \forall\ t\in[0,T],
\end{array}
\right.
\end{equation}
where $\tilde \sigma_{ij} = \sum_{k=1}^n \sigma_{ik}(X) \sigma_{jk}(X)$.

Essentially, the partial differential equation (PDE) (\ref{mainHJI}) (a) is the Hamilton-Jacobi-Isaacs equation which returns the value function $v(x,y,m,t)$ once we fix the distribution $m(x,y,t)$; This PDE has to be solved backwards with boundary conditions at final time $T$, represented by the last line in \ref{mainHJI} (a). In \ref{mainHJI} (b) we have the optimal closed-loop control $u^*(x,t)$ and worst-case disturbance $w^*(x,t)$ as min-maximizers of the Hamiltonian function in the RHS. The PDE (\ref{mainFPK}) represents the transport equation of the measure $m$ immersed in a vector field $AX + Bu +C+D w$; It returns the distribution $m(x,y,t)$ once fixed both $u^*(x,t)$ and $w^*(x,t)$  and consequently the vector field $AX + Bu^* +C+D w^*$. Such a PDE has to be solved forwards with boundary condition at the initial time (see the fourth line) of \ref{mainFPK}. Finally, once given $m(x,y,t)$ from (c) and entered into the running cost $c(x,y,m,u)$ in (a), we obtain the error
\begin{equation}
\label{mainFPK1}
\left\{
\begin{array}{ll}
m_{on}(t) := \int_{x_{on}}^{x_{off}} \int_{[0,1]} y m(x,y,t)dxdy \ \ \forall\ t\in[0,T],\\
e(t)=m_{on}(t) - \overline m_{on}.
\end{array}
\right.
\end{equation}

Note that 
$$  \bar X(t)=\left[\begin{array}{ll}
\bar x(t) \\
\bar y(t)\end{array}\right] = \left[\begin{array}{ll}
\bar x(t) \\
m_{on}\end{array}\right] =\left[\begin{array}{ll}
\int_{x_{on}}^{x_{off}} \int_{[0,1]} x m(x,y,t)dxdy  \\
\int_{x_{on}}^{x_{off}} \int_{[0,1]} y m(x,y,t)dxdy\end{array}\right], 
$$
and therefore, henceforth we can refer to as mean-field equilibrium solutions any  pair
$(v(X,t),\bar X(t))$ which is solution of (\ref{mainHJI})-(\ref{mainFPK}).

\section{Main results}\label{mainresults}

This paper contributes in three directions with respect to the TCLs application introduced earlier. First, it analyzes and computes the  mean-field equilibrium for the deterministic mean-field game and proves that under certain conditions the microscopic dynamics is asymptotically stable.
 We repeat the analysis for the stochastic case, assuming that the microscopic dynamics is uncertain. 
  Even for this case, a mean-field equilibrium is computed, and stochastic stability is studied. We distinguish two cases. On the one hand, we consider a stochastic disturbance  which fades to zero the closer the state is to zero. The Brownian motion coefficients are linear in the state and the resulting dynamics is also known as   geometric Brownian motion. 
     On the other hand, 
we take the stochastic disturbance being independent on the state. The Brownian motion coefficients are constant and the resulting dynamics mirrors the Langevin equation. In both cases we prove stochastic stability of second-moment for the stochastic process at hand. This section ends with a detailed analysis of robustness properties. The microscopic dynamics is now subject to an addition exogenous input, the disturbance, with bounded resource or energy. Even for this case, we study the mean-field equilibrium and  investigate on condition that guarantee stability.

\subsection{Mean-field equilibrium and stability}

In this section we  establish an explicit solution in terms of mean-field equilibrium  for the  deterministic case and study stability of the microscopic dynamics. This case is obtained by fixing to zero the coefficients of both stochastic and adversarial disturbance. 

The  linear quadratic problem we wish to solve is then:
\begin{equation}\label{pq}
\begin{array}{cll}
\displaystyle
\inf_{\{u_t\}_{t}}  \int_0^T \left[\frac{1}{2} \Big(X(t)^T  Q X(t) +  u(t)^T R u(t)^T\Big)+ L^T X(t) \right] dt\\ \\
\dot X(t) = A X(t) +  B u(t) + C   \mbox{ in } \mathcal S.
\end{array}
\end{equation}
%

The next result shows that the problem reduces to solving three matrix equations. 

\begin{theorem}\label{thm1} \textbf{(Mean-field equilibrium)}
Let $D,\Sigma=0$ in the game (\ref{mainHJI})-(\ref{mainFPK}). A mean-field equilibrium for  (\ref{mainHJI})-(\ref{mainFPK}) is given by 
\begin{equation}\label{mf}\displaystyle
\left\{
\begin{array}{l}v(X,t)= \frac{1}{2}X^TP(t) X+ \Psi(t)^T X + \chi(t),\\ \\
\dot{\bar X}(t) = [A(x) - B R^{-1} B^T P ] \bar X(t) - B R^{-1} B^T \bar \Psi(t)  + C,
\end{array}\right.
\end{equation}
where 
\begin{equation}\label{cb4state}
\left\{
\begin{array}{l}
\dot P  + PA(x) + A(x)^T P - P BR^{-1} B^T P + Q=0
\ \mbox{in }  [0,T[,\ P(T) =\phi, \\ \\
\dot \Psi + A(x)^T \Psi + P C -  P BR^{-1} B^T \Psi + L =0 \ \mbox{in } [0,T[,\ \Psi(T) =0, \\ \\
\dot \chi  + \Psi^T C - \frac{1}{2} \Psi^T B R^{-1} B^T \Psi =0  \ \mbox{in } [0,T[,\ \chi(T) =0,
\end{array}\right.
\end{equation}
and $\bar \Psi(t)=\int_{x_{on}}^{x_{off}} \int_{[0,1]} \Psi(t) m(x,y,t)dxdy$. Furthermore, the mean-field equilibrium strategies are given by
\begin{equation}\label{mfes} 
u^*(X,t) = - R^{-1} B^T [P X + \Psi].\end{equation}
\end{theorem}

\textbf{Proof.} Given in the appendix. \qed 

Let us note that by substituting the mean-field equilibrium strategies $u^* = - R^{-1} B^T [P X + \Psi]$  given in (\ref{mfes}) in the open-loop microscopic dynamics $\dot X(t) = A X(t) +  B u(t) + C$ as defined in (\ref{pq}), the closed-loop microscopic dynamics  is
\begin{equation}\label{asdyn}\dot X(t) = [A(x) - B R^{-1} B^T P ] X(t) - B R^{-1} B^T \Psi(x,e,t)  + C.\end{equation}

Now,  let $\mathcal X$ be the set of equilibrium points for (\ref{asdyn}), namely, the set of $X$ such that 
$$\mathcal X=\{(X,e) \in \mathbb R^2 \times \mathbb R|\, [A(x) - B R^{-1} B^T P ] X(t) - B R^{-1} B^T \Psi(x,e,t)  + C=0\},$$
and let   $V(X(t)) = dist(X(t),\mathcal X)$.
The next result establishes a condition under which  the above dynamics converges asymptotically to the set of equilibrium points. 

\begin{corollary}\label{cor1} \textbf{(Asymptotic stability)}
If  it holds
\begin{equation}
\begin{array}{ll}
\partial_X V(X,t)^T \Big([A - B R^{-1} B^T  P] X(t)  - B R^{-1} B^T \Psi^*(x(t),e(t))  + C\Big) \\
 \qquad \qquad \qquad  < - \|X(t) - \Pi_{\mathcal X}(X(t))\|^2 
 \end{array}
 \end{equation}
then dynamics (\ref{asdyn}) is asymptotically stable, namely, $\lim_{t \rightarrow \infty} dist(X(t),\mathcal X)=0$.
\end{corollary}

\textbf{Proof.} Given in the appendix. \qed 

\subsection{Stochastic case}
In this section we study the case where the dynamics is given by a stochastic differential equation driven by a Brownian motion.  In other words, the model is uncertain and the uncertainty is modeled as a stochastic disturbance. 

The problem at hand is then:
\begin{equation}\label{pq1}
\begin{array}{cll}
\displaystyle
\inf_{\{u_t\}_{t}}  \mathbb E \int_0^T \left[\frac{1}{2} \Big(X(t)^T  Q X(t) +  u(t)^T R u(t)^T\Big)+ L^T X(t) \right] dt\\ \\
d X(t) = (A X(t) +  B u(t) + C) dt + \Sigma d\mathcal B_t,   
\end{array}
\end{equation}
where all matrices are as in (\ref{fd1}) and 
\begin{equation}\label{fd2}\nonumber
\begin{array}{cccc}
  \Sigma = \left[\begin{array}{ccc}
\sigma_{11}(x) & 0\\
0 &  \sigma_{22}(y) \end{array}\right].
\end{array}
\end{equation}

This section investigates on the solution of the HJI equation under the assumption that the time evolution of the common state is given. We show that the problem reduces to solving three matrix equations. 
To see this, by isolating the HJI part of (\ref{mainHJI}) for fixed $m_t$, for $t \in [0,T]$, we have

\begin{equation}\label{stochmfg}\nonumber
\left\{
\begin{array}{lll}
-\partial_t v(X,t) - \sup_u \Big\{ - \partial_X v(X,t)^T (AX + Bu +C) - \frac{1}{2} \Big(X^T  Q X \\
 \qquad \qquad-  u^T R u\Big)- L^T X \Big\}+ \frac{1}{2} ( \sigma_{11}(x)^2  \partial_{xx} v(X,t) \\  \qquad \qquad  \qquad 
 +\sigma_{22}(y)^2  \partial_{yy} v(X,t))=0,   \mbox{ in } \mathcal S \times [0,T[,\\
v(X,T)=g(x)   \mbox{ in } \mathcal S,\\\\
\displaystyle 
u^*(x,t)=-r^{-1} B^T \partial_y v(X,t).
\end{array}\right.
\end{equation}

\noindent

Let us consider the following value function 
$$v(X,t)= \frac{1}{2}X^TP(t) X+ \Psi(t)^T X + \chi(t),$$ 

and 
$$u^* = - R^{-1} B^T [P X + \Psi],$$

so that (\ref{cb2}) can be rewritten as

\begin{equation}\label{cb333}\displaystyle
\left\{
\begin{array}{r}
\frac{1}{2} X^T\dot P(t) X + \dot \Psi(t)X + \dot \chi(t) +(P(t) X+ \Psi(t))^T \Big[ -BR^{-1} B^T \Big]  
(P(t) x+ \Psi(t))  \\
+ (P(t) x+ \Psi(t))^T (AX +C) 
+ \frac{1}{2} \Big(X(t)^T  Q X(t) 
 +  u(t)^T R u(t)^T\Big) \\ 
  + L^T X(t) +\frac{1}{2} ( \sigma_{11}(x)^2 P_{11}(t)+\sigma_{22}(y)^2 P_{22}(t))=0 
\   \mbox{ in } \mathcal S\times[0,T[,
\\ \\
P(T) =\phi, \quad \Psi(T) =0, \quad \chi(T)=0.
\end{array}\right.
\end{equation}

%




The boundary conditions are obtained by imposing that $$v(X,T)= \frac{1}{2} X^T P(T) X+ \Psi(T) X + \chi(T)=\frac{1}{2}X^T \phi X.$$

\subsubsection{Case I: state dependent variance}
The first case we consider involves coefficients for the Brownian motion linear in the state, namely 
\begin{equation}\label{case1}\Sigma(X) = \left[\begin{array}{ccc}
\hat \sigma_{11}x & 0\\
0 & \hat  \sigma_{22}y \end{array}\right]\end{equation}

\begin{theorem}\label{thm2} \textbf{(stochastic mean-field equilibrium: case I)}
A mean-field equilibrium for the game (\ref{eq:meanfieldprel2})  with $\Sigma(X)$ as in (\ref{case1}) is given by 
\begin{equation}\label{mf}\displaystyle
\left\{
\begin{array}{l}v(X,t)= \frac{1}{2}X^TP(t) X+ \Psi(t)^T X + \chi(t),\\ \\
\dot{\bar X}(t) = [A - B R^{-1} B^T P ] \bar X(t) - B R^{-1} B^T \Psi^*(\bar X(t))  + C,
\end{array}\right.
\end{equation}
where 
\begin{equation}\label{cb4statestoch}
\left\{
\begin{array}{l}
\dot P(t)  + P(t)A + A^T P - P BR^{-1} B^T P + Q+ \tilde P=0
\ \mbox{in }  [0,T[,\ P(T) =\phi, \\ \\
\dot \Psi(t) + A^T \Psi + P C -  P BR^{-1} B^T \Psi + L =0 \ \mbox{in } [0,T[,\ \Psi(T) =0, \\ \\
\dot \chi(t)  + \Psi(t)^T C - \frac{1}{2} \Psi^T B R^{-1} B^T \Psi =0  \ \mbox{in } [0,T[,\ \chi(T) =0,
\end{array}\right.
\end{equation}
and 
\begin{equation}\label{diag}
\tilde P = Diag( (\hat \sigma_{ii}^2  P_{ii})_{i =1,2} )=
\left[\begin{array}{cc} \hat \sigma_{11}^2  P_{11} & 0 \\
0 & \hat \sigma_{22}^2 P_{22} \end{array}\right].\end{equation}
Furthermore, the mean-field equilibrium strategy is 
\begin{equation}
\label{mfsestoch}u^* = - R^{-1} B^T [P X + \Psi]
\end{equation}
\end{theorem}

\textbf{Proof.} Given in the appendix. \qed

Based on the above result, let us now substitute the expression of the mean-field equilibrium strategy $u^* = - R^{-1} B^T [P X + \Psi]$ as in (\ref{mfsestoch}) in the open-loop microscopic dynamics $d X(t) = (A X(t) +  B u(t) + C) dt + \Sigma d\mathcal B(t)$ given in (\ref{pq1}) so to obtain
the closed-loop microscopic dynamics
\begin{equation}\label{asdyn1111} 
d X(t) =\Big[ (A(x) - B R^{-1} B^T P ) X(t) - B R^{-1} B^T \Psi^*(x(t),e(t))  + C \Big] dt+ \Sigma d\mathcal B(t) \end{equation}

Now,  let $\mathcal X$ be the set of equilibrium points for (\ref{asdyn1111}), namely, the set of $X$ such that 
$$\mathcal X=\{(X,e) \in \mathbb R^2 \times \mathbb R|\, (A(x) - B R^{-1} B^T P ) X(t) - B R^{-1} B^T \Psi^*(x,e)  + C=0\},$$
and let   $V(X(t)) = dist(X(t),\mathcal X)$.
The next result establishes a condition under which  the above dynamics converges asymptotically to the set of equilibrium points. 

\begin{corollary}\label{cor2} \textbf{(2nd moment boundedness)}
Let a compact set $\mathcal M\subset  \mathbb R^2$ be given. 
Suppose that for all $X \not \in \mathcal M$ 
\begin{equation}\label{contract}
\begin{array}{lll} \partial_X V(X,t)^T \Big([A - B R^{-1} B^T  P] X(t) 
- B R^{-1} B^T \Psi^*(x(t),e(t))  + C\Big) 
  \\   \qquad \qquad  \qquad  < -\frac{1}{2} ( \sigma_{11}^2(x)  \partial_{xx} V(X,t)
 +  \sigma_{22}^2(x)  \partial_{yy} V(X,t) )
 \end{array}
 \end{equation}



then  dynamics (\ref{asdyn1111}) is a stochastic process with 2nd moment bounded.
\end{corollary}
\textbf{Proof.} Given in the appendix. \qed

\subsubsection{Case II: state independent variance and Langevin equation}
The second case we consider involves coefficients for the Brownian motion which are constant, namely
\begin{equation}\label{case2}\Sigma = \left[\begin{array}{ccc}
\hat \sigma_{11} & 0\\
0 & \hat  \sigma_{22} \end{array}\right] \end{equation}

\begin{theorem} \label{thm3} \textbf{(stochastic mean-field equilibrium: case II)}

Let $\Sigma$ be as in (\ref{case2}). 
A mean-field equilibrium for the game (\ref{mainHJI})-(\ref{mainFPK})  is given by 
%
\begin{equation}\label{mf}\displaystyle
\left\{
\begin{array}{l}v(X,t)= \frac{1}{2}X^TP(t) X+ \Psi(t)^T X + \chi(t),\\ \\
\dot{\bar X}(t) = [A - B R^{-1} B^T P ] \bar X(t) - B R^{-1} B^T \Psi^*(\bar X(t))  + C
\end{array}\right.
\end{equation}
where 
\begin{equation}\label{cb4stoch2}
\left\{
\begin{array}{l}
\dot P(t)  + P(t)A + A^T P - P BR^{-1} B^T P + Q=0
\ \mbox{in }  [0,T[,\ P(T) =\phi, \\ \\
\dot \Psi(t) + A^T \Psi + P C -  P BR^{-1} B^T \Psi + L =0 \ \mbox{in } [0,T[,\ \Psi(T) =0, \\ \\
\dot \chi(t)  + \Psi(t)^T C - \frac{1}{2} \Psi^T B R^{-1} B^T \Psi + \tilde P =0  \ \mbox{in } [0,T[,\ \chi(T) =0,
\end{array}\right.
\end{equation}
and 
\begin{equation}\label{diagstoch2}
\tilde P = 
\left[\begin{array}{cc} \hat \sigma_{11}^2& 0 \\
0 & \hat \sigma_{22}^2  \end{array}\right].\end{equation}
Furthermore, the mean-field equilibrium strategies are given by
\begin{equation}\label{mfes1} 
u^*(X,t) = - R^{-1} B^T [P X + \Psi].\end{equation}
\end{theorem}

\textbf{Proof.} Given in the appendix. \qed
%
%
%

Based on the above result, let us now substitute the expression of the mean-field equilibrium strategy $u^* = - R^{-1} B^T [P X + \Psi]$ as in (\ref{mfes1}) in the open-loop microscopic dynamics $d X(t) = (A X(t) +  B u(t) + C) dt + \Sigma d\mathcal B(t)$ given in (\ref{pq1}) so to obtain
the closed-loop microscopic dynamics
\begin{equation}\label{asdyn11} 
d X(t) =\Big[ (A(x) - B R^{-1} B^T P ) X(t) - B R^{-1} B^T \Psi^*(x(t),e(t))  + C \Big] dt+ \Sigma d\mathcal B(t) \end{equation}

%
Now,  let $\mathcal X$ be the set of equilibrium points for (\ref{asdyn}), namely, the set of $X$ such that 
$$\mathcal X=\{(X,e) \in \mathbb R^2 \times \mathbb R|\, (A(x) - B R^{-1} B^T P ) X(t) - B R^{-1} B^T \Psi^*(x,e)  + C=0\},$$
and let   $V(X(t)) = dist(X(t),\mathcal X)$.
The next result establishes a condition under which  the above dynamics converges asymptotically to the set of equilibrium points. 

\begin{corollary} \label{cor3} \textbf{(2nd moment boundedness)}
Let a compact set $\mathcal M\subset  \mathbb R^2$ be given. 
Suppose that for all $X \not \in \mathcal M$
\begin{equation}\label{contract1}
\begin{array}{ll}
\partial_X V(X,t)^T \Big([A - B R^{-1} B^T  P] X(t)  \\ \qquad \qquad - B R^{-1} B^T \Psi^*(x(t),e(t))  + C\Big)< 
-\frac{1}{2} ( \hat \sigma_{11}^2  \partial_{xx} V(X,t)
 + \hat \sigma_{22}^2  \partial_{yy} V(X,t) )
 \end{array}\end{equation}

then  dynamics (\ref{asdyn11}) is a stochastic process with 2nd moment bounded.
\end{corollary}
\textbf{Proof.} Given in the appendix. \qed

\subsection{Model miss-specification}
This section deals with model miss-specification, this being represented by an additional exogenous and adversarial disturbance. The disturbance is supposed to be of bounded energy. 
Thus,  the  linear quadratic problem we wish to solve is:
\begin{equation}\label{pq4}
\begin{array}{cll}
\displaystyle
\inf_{\{u_t\}_{t}}  \mathbb E \int_0^T \left[\frac{1}{2} \Big(X(t)^T  Q X(t) +  u(t)^T R u(t) - \gamma^2 w(t)^T  w(t) \Big)+ L^T X(t) \right] dt\\ \\
\dot X(t) = A X(t) +  B u(t) + C + Dw(t) \mbox{ in } \mathcal S.   
\end{array}
\end{equation}
%

This section investigates on the solution of the HJI equation under the assumption that the time evolution of the common state is given. We show that the problem reduces to solving three matrix equations. 
To see this, by isolating the HJI part of (\ref{mainHJI}) for fixed $m_t$, for $t \in [0,T]$, we have
%

\begin{theorem}\label{thm4} \textbf{(worst-case mean-field equilibrium)}
A mean-field equilibrium for  (\ref{mainHJI})-(\ref{mainFPK}) is given by 
\begin{equation}\label{mf}\displaystyle
\left\{
\begin{array}{l}v(X,t)= \frac{1}{2}X^TP(t) X+ \Psi(t)^T X + \chi(t),\\ \\
\dot{\bar X}(t) = [A - B R^{-1} B^T P ] \bar X(t) - B R^{-1} B^T \Psi^*(\bar X(t))  + C,
\end{array}\right.
\end{equation}
where 
\begin{equation}\label{cbwc}
\left\{
\begin{array}{l}
\dot P(t)  + P(t)A + A^T P +P (- B R^{-1} B^T  + \frac{1}{\gamma^2} DD^T) P + Q=0
\ \mbox{in }  [0,T[,\ P(T) =\phi, \\ \\
\dot \Psi(t) + A^T \Psi + P C +(- B R^{-1} B^T  + \frac{1}{\gamma^2} DD^T) \Psi + L =0 \ \mbox{in } [0,T[,\ \Psi(T) =0, \\ \\
\dot \chi(t)  + \Psi(t)^T C + \frac{1}{2} \Psi^T (- B R^{-1} B^T  + \frac{1}{\gamma^2} DD^T) \Psi =0  \ \mbox{in } [0,T[,\ \chi(T) =0,
\end{array}\right.
\end{equation}
Furthermore, the mean-field equilibrium control and disturbance are
\begin{equation}\label{mfes4}
\begin{array}{l}
u^* = - R^{-1} B^T [P X + \Psi]\\
w^* = \frac{1}{\gamma^2} D^T [P X + \Psi].
\end{array}\end{equation}
\end{theorem}

\textbf{Proof.} Given in the appendix. \qed

Let us note that by substituting the mean-field equilibrium strategies $u^* = - R^{-1} B^T [P X + \Psi]$  and $w^* = \frac{1}{\gamma^2} D^T [P X + \Psi]$ as given in (\ref{mfes4}) in the open-loop microscopic dynamics $\dot X(t) = A X(t) +  B u(t) + C + Dw$ as defined in (\ref{pq4}), the closed-loop microscopic dynamics  is
\begin{equation}\label{wcdyn}
\dot X(t) = [A(x) +(- B R^{-1} B^T  + \frac{1}{\gamma^2} DD^T)P] X(t) +(- B R^{-1} B^T +  \frac{1}{\gamma^2} DD^T)\Psi^*(x(t),e(t))  + C \end{equation}

Now,  let $\mathcal X$ be the set of equilibrium points for (\ref{asdyn}), namely, the set of $X$ such that 
$$\mathcal X=\{(X,e) \in \mathbb R^2 \times \mathbb R|\, [A(x) +(- B R^{-1} B^T + \frac{1}{\gamma^2} DD^T) P ] X(t) - B R^{-1} B^T \Psi(x,e,t)  + C=0\},$$
and let   $V(X(t)) = dist(X(t),\mathcal X)$.
The next result establishes a condition under which  the above dynamics converges asymptotically to the set of equilibrium points.

%
%

\begin{corollary} \label{cor4} \textbf{(worst-case stability)}
If it holds
\begin{equation}
\begin{array}{ll}
\partial_X V(X,t)^T \Big([A +(- B R^{-1} B^T  + \frac{1}{\gamma^2} DD^T)P] X(t) +(- B R^{-1} B^T +  \frac{1}{\gamma^2} DD^T)\\
 \qquad \qquad \qquad \cdot \Psi^*(x(t),e(t))  + C\Big) 
  < - \|X(t) - \Pi_{\mathcal X}(X(t))\|^2 
 \end{array}
 \end{equation}
then dynamics (\ref{wcdyn}) is asymptotically stable, namely, $\lim_{t \rightarrow \infty} dist(X(t),\mathcal X)=0$.

\end{corollary}
\textbf{Proof.} Given in the appendix. \qed

\section{Numerical studies}\label{sec:num}
In this section a system consisting of  $n=10^2$ indistinguishable TCLs. All simulations are carried out with MATLAB on an Intel(R) Core(TM)2 Duo, 
CPU P8400 at 2.27 GHz and a 3GB of RAM.  The number of iterations is $T=30$.
We consider a discrete time version of (\ref{pq}) 
 \begin{equation}\label{pqdisc}X(t+dt) = X(t)+(A(x(t)) X(t) +  B u(t) + C)dt. \end{equation}

The parameter  are as shown in Table \ref{table:param} and in particular  the step size $dt=0.1$, the cooling and heating rates are $\alpha=\beta=1$, the lowest and highest temperatures are $x_{on}=-10$, 
and $x_{off}=10$, respectively, the penalty coefficients are $r_{on}=r_{off}=1$, and $q=1$, and the initial distribution is normal with zero mean and standard deviation $std(m(0))=1$.
\begin{table}
\begin{center}
   \begin{tabular}{|c|c|c|c|c|c|c|c|c|}
  \hline
      $\alpha$ & $\beta$ &  $x_{on}$ & $x_{on}$  & $r_{on},\ r_{off}$ & $q$ & $std(m_0)$ & $\bar m_0$ \\ \hline \hline
        $1$  &   $1$    & $-10$ & $10$  & $10$     & $1$ & $1$   & $0$  \\ \hline \end{tabular}\\
  \end{center}
  \caption{Simulation parameters }\label{table:param}\end{table}

 The numerical results are obtained using the algorithm in Table \ref{fig:algorithm} for a discretized set of states. 


\begin{table}\normalsize
\begin{center}
\begin{tabular}{p{12.5cm}}\\
\toprule
\textbf{Input:} Set of parameters as in Table \ref{table:param}.  \\
\textbf{Output:} TCLs' states $X(t)$\\
$\quad 1: $ \textbf{Initialize.} Generate $X(0)$ given $\bar m_0$ and $\mbox{std}(m_0)$\\
$\quad 2: $ \textbf{for} time $iter=0,1,\ldots,T-1$ \textbf{do}\\
$\quad 3: \quad$ \textbf{if} $iter>0$, \textbf{then} compute $m_t$, $\bar m_t$, and $\mbox{std}(m_t)$\\
$\quad 4: \quad$ \textbf{end if} \\
$\quad 5: \quad$ \textbf{for} player $i=1,\ldots,n$ \textbf{do}\\
$\quad 6: \quad$ Set $t=iter \cdot dt$ and compute control $\tilde u(t)$ using current $\bar m(t)$  \\
$\quad 7: \quad\quad$ compute new state  $X(t+dt)$ by executing (\ref{pqdisc})\\
$\quad 8: \quad$ \textbf{end for}\\
$\quad 9: $ \textbf{end for} \\
$\quad 10: $ \textbf{STOP}\\
\bottomrule
\end{tabular}
\end{center}\caption{Simulation algorithm} \label{fig:algorithm}
\end{table}


The optimal control is taken as
$$u^* = - R^{-1} B^T [P X + \Psi] $$
where $P$ is obtained from running the  MATLAB command \texttt{[P]=care(A,B,Q,R)}, which receives the matrices as input and returns the solution $P$ to the algebraic  Riccati equation.
Under the assumption  $B R^{-1} B^T  \Psi \approx C$
the resulting closed-loop dynamics is given by
$$X(t+dt) = X(t)+ [A - B R^{-1} B^T P ] X(t) dt.$$

\begin{figure} [htb]
\centering
\includegraphics[width=\columnwidth]{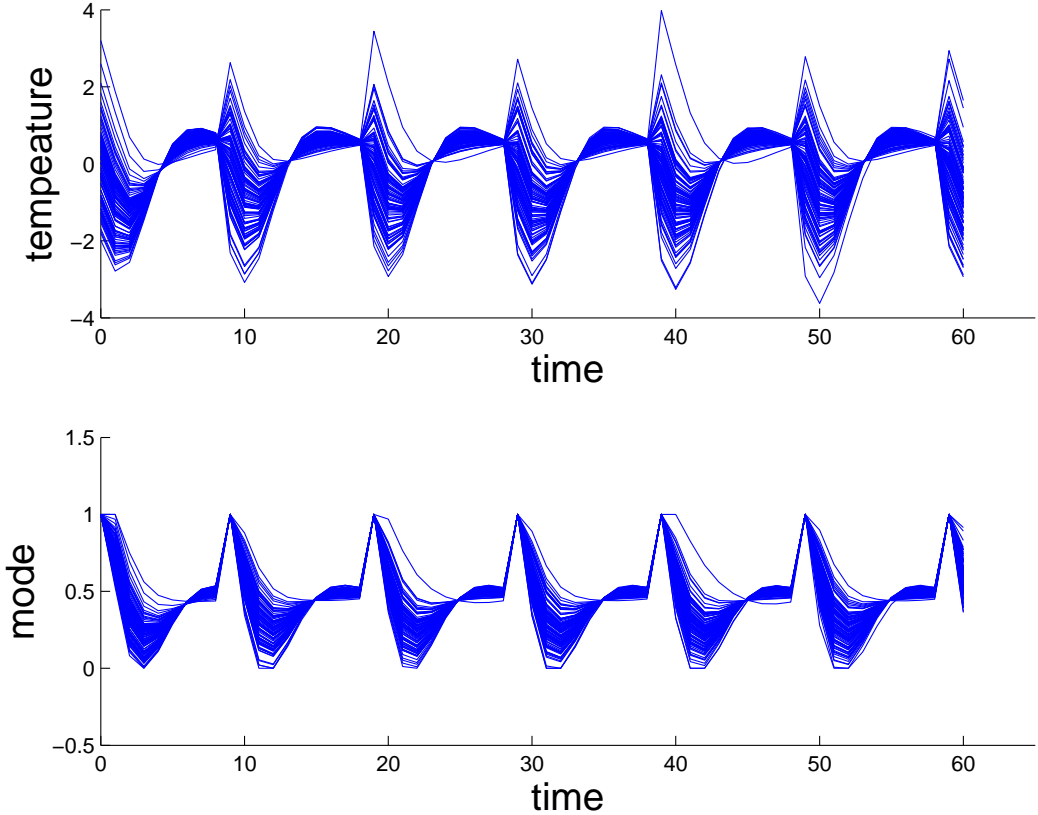}
\caption{Time plot  of the state of each TCL, namely temperature $x(t)$ (top row) and mode $y(t)$ (bottom row).}
\label{fig:Fig1}
\end{figure}

Figure \ref{fig:Fig1} displays the time plot  of the state of each TCL, namely its temperature $x(t)$ (top row) and mode $y(t)$ (bottom row). In contrast with what we observed in Fig. \ref{fig:oscillations}, the TCLs show a stable behavior. The simulation is carried out assuming that any 10 seconds the states are subject to an impulse. The TCLs react to the impulse very fast  and converge to the equilibrium point before a new impulse is activated, as clear visually in the plot.

\begin{figure} [htb]
\centering
\includegraphics[width=\columnwidth]{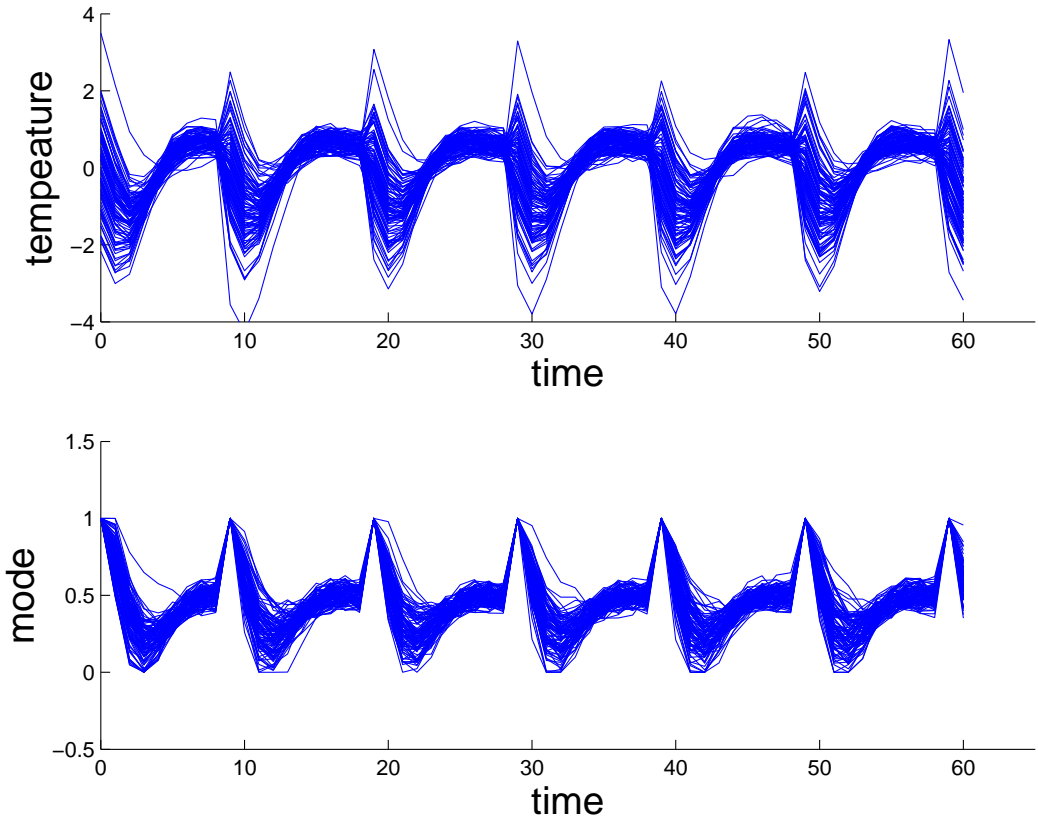}
\caption{Time plot  of the state of each TCL, namely temperature $x(t)$ (top row) and mode $y(t)$ (bottom row).}
\label{fig:Fig2}
\end{figure}

We repeat the simulation for the two stochastic cases discussed earlier. The stochastic version of the dynamics appears now as
$$X(t+dt) = X(t)+(A(x(t)) X(t) +  B u(t) + C + \Sigma W(t))dt $$ or
for the first case, and 
 $$X(t+dt) = X(t)+(A(x(t)) X(t) +  B u(t) + C + \Sigma(x) W(t))dt $$
 for the second case. Here  $W(t)$ is a random walk. 
The corresponding closed-loop dynamics are then  $$X(t+dt) = X(t)+ [A - B R^{-1} B^T P ] X(t) dt+ \Sigma W(t) dt $$ and
 $$X(t+dt) = X(t)+ [A - B R^{-1} B^T P ] X(t) dt+ \Sigma(x) W(t))dt $$ respectively. 
Figure \ref{fig:Fig2} displays the time plot  of the state of each TCL, namely its temperature $x(t)$ (top row) and mode $y(t)$ (bottom row) in the first case. Even in this case, differently from what observed sin Fig. \ref{fig:oscillations}, the TCLs  react to the impulse  and converge to the equilibrium point before a new impulse is activated. The effects of the Brownian motion is the one of enlarging the domain of attraction.

\begin{figure} [htb]
\centering
\includegraphics[width=\columnwidth]{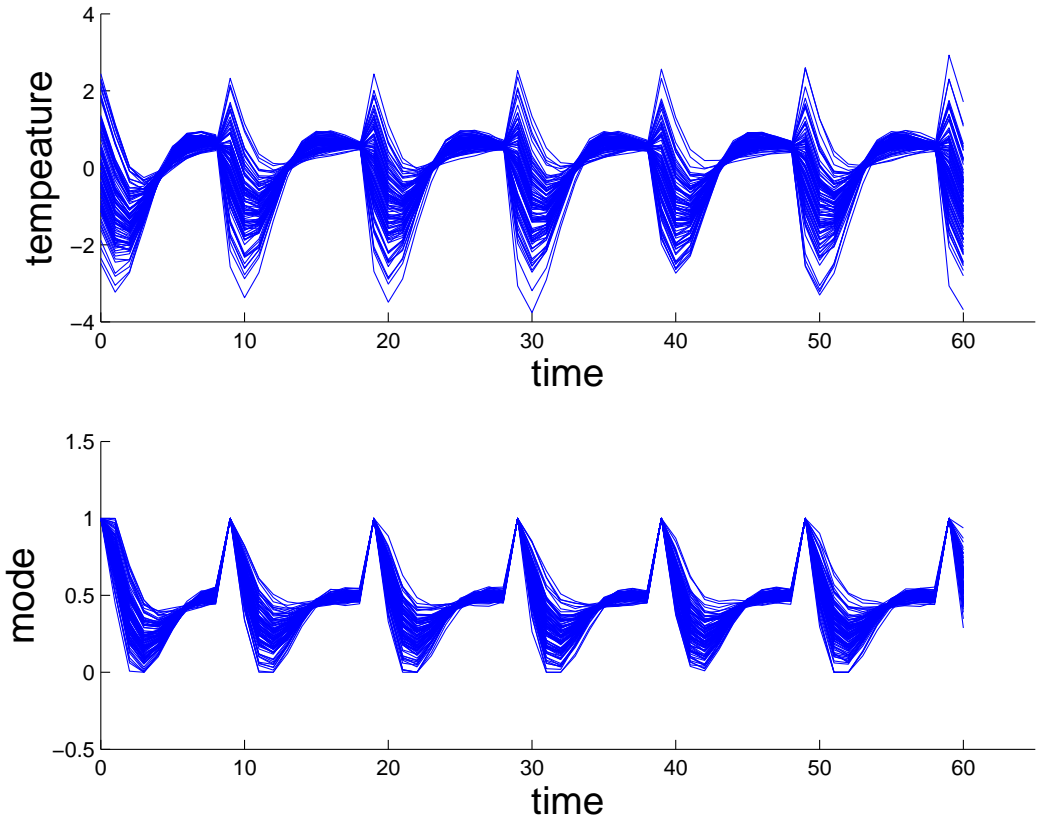}
\caption{Time plot  of the state of each TCL, namely temperature $x(t)$ (top row) and mode $y(t)$ (bottom row).}
\label{fig:Fig3}
\end{figure}

The experiment is repeated in Figure \ref{fig:Fig3} for the geometric Brownian motion. As in the previous cases the figure  displays the time plot  of the state of each TCL, namely its temperature $x(t)$ (top row) and mode $y(t)$ (bottom row) in the first case. As the Brownian motion is not weighted by the state (in modulus), its effects are attenuated and the plot is more similar to the one in Fig. \ref{fig:Fig1}.

%

Note that except for the Langevin-type dynamics, in  the remainder  two cases the TCLs states are driven to  zero. For the Langevin-type dynamics the state is confined within a neighborhood of zero. 



\section{Discussion} \label{discussion}
With regards to the problem at hand, the topic of dynamic response management has sparked the attention of scientists from different disciplines. This is witnessed by the rapid growing of publications in journals of different research areas, from differential game theory \cite{BagBau13,CPTD12,PCGL14}, to control and optimization \cite{AK12,CH11,MCH13,MKC13,RDM12}, to computer science \cite{bcSGEFA14}. One reason for this is that dynamic response management intersects research programs in \emph{smart buildings} and \emph{smart cities}. The problem is relevant due to an ever-increasing size of network systems and the consequent impossibility of centralizing the management of the whole system.

Fully aware of the importance of the topic, let us discuss the  relevance of the results of this paper. First, it must be said that the game-theoretic approach presented here is a natural way to deal with larges scale, complex and distributed systems where no central planner may be capable of processing all information data and in order to control the whole system online. One way to deal with this issue, and  which is the main idea of dynamic demand management, aims at assigning part of the regulation burden to the consumers by using frequency responsive appliances. In other words, each appliance regulates automatically and in a decentralized fashion its power demand based on the mains frequency.

In this respect, the provided model builds upon the strategic interaction among the electrical appliances. Note that here we look at the problem in more general terms and talk about electrical appliances rather than TCLs. The model suits the case where where the latter are numerous and indistinguishable. Indistinguishable means that any appliance in the same condition will react at the same way. We wish to highlight  that indistinguishability is not a limitation, as in the case of heterogeneity of the electrical appliances, more complex multi-population models may be derived based on the same modeling approach used here. 

The results provided in this paper shed light on the existence of mean-field equilibrium solutions.  By this we mean strategies based on the current and  forecasted demand, which are proven to attenuate oscillations of mains frequency. A first feature of the model at hand is that the considered strategies are stochastic. This means that the TCL sets a probability with which to switch $on$ or $off$. Stochastic linear strategies are designed as closed-loop feedback strategies on current state, temperature and switching mode. Such strategies are computed over a finite horizon and therefore are based  on forecasted demand. From another angle, we may say that mean-field equilibrium strategies represent the asymptotic limit of Nash equilibrium strategies, and as such they are the best-response strategies of even single player, for fixed behavior of other players. The proven stability of the microscopic dynamics confirms the asymptotic convergence of the TCLs's states to an equilibrium point, this being expressed  in terms of temperature and switching  mode. The several cases studied in the paper have shown that this holds true in the cases of both perfect and imperfect modeling. This is a clear evidence of a certain degree of robustness characterizing  the proposed strategies. In the case of imperfect modeling, model mis-specifications is considered both in a stochastic and deterministic worst-case scenario. Assuming imperfect models both with a stochastic or worst-case deterministic disturbance acting on the state dynamics, conditions for convergence of the microscopic dynamics are provided.

 \section{Concluding remarks}\label{conclusion}

%
%

We have illustrated robust mean-field games as a paradigm for crowd-averse systems. We have discussed these systems in the context of stock market, production engineering, and dynamic demand management in power systems. As main contributions we first have formulated the problem as a robust mean-field game; second, we have presented  a new approximation method  based on the extension of the state space; third we have discussed a relaxation method to minimize the approximation error. Further results are obtained for a scalar microscopic dynamics, for which we have established performance bounds, and analyzed stochastic stability of both the microscopic and the macroscopic dynamics. 
We can extend our study in at least three directions. These include i)  the extension  of the approximation method  to more general cost functionals, ii) the study of the case with ``local'' mean-field interactions rather than ``global'' as in the current scenario, and iii) the analysis of crowd-seeking  scenarios in contrast to the crowd-averse cases analyzed in this paper.

\section*{Appendix}

\subsection*{Proof of Theorem \ref{thm1}}

Let us start by isolating the HJI part of (\ref{mainHJI}). For fixed $m_t$ and for $t \in [0,T]$, we have
\begin{equation}\label{cb2}
\left\{
\begin{array}{l}
\displaystyle
-\partial_t v(x,y,t)- \Big\{ y\Big[ -\alpha (x-x_{on})\Big] + (1-y)\Big[ -\beta (x-x_{off})\Big] \Big\} \partial_x v(x,y,t)  \\+ \sup_{u\in \mathbb R} \Big\{-B u \,\partial_y v(x,y,t) - \frac{1}{2} q x^2  + \frac{1}{2} u^T r u + y (S e  + W)\Big\}=0 \\
   \mbox{ in } \mathcal S \times]0,T],\\
\displaystyle
v(x,y,T)=\Psi(x)\   \mbox{ in } \mathcal S,\\\\
\displaystyle 
u^*(x,t)=-r^{-1} B^T \partial_y v(x,y,t)
\end{array}\right.
\end{equation}
which in a more compact form can be rewritten as
\begin{equation}\label{1111expl}\nonumber
\left\{
\begin{array}{lll}
-\partial_t v(X,t) - \sup_u \Big\{  \partial_X v(X,t)^T (AX + Bu +C) + \frac{1}{2} \Big(X^T  Q X \\
 \qquad \qquad +  u^T R u^T\Big)+ L^T X \Big\}=0,   \mbox{ in } \mathcal S \times [0,T[,\\ 
v(X,T)=g(x)   \mbox{ in } \mathcal S,\\\\
\displaystyle 
u^*(x,t)=-r^{-1} B^T \partial_y v(X,t).
\end{array}\right.
\end{equation}

\noindent

Let us consider the following value function 
$$v(X,t)= \frac{1}{2}X^TP(t) X+ \Psi(t)^T X + \chi(t),$$ 

and the corresponding optimal closed-loop state feedback strategy
$$u^* = - R^{-1} B^T [P X + \Psi].$$

Then (\ref{cb2}) can be rewritten as
\begin{equation}\label{cb33}\displaystyle
\left\{
\begin{array}{r}
\frac{1}{2} X^T\dot P(t) X + \dot \Psi(t)X + \dot \chi(t) +(P(t) X+ \Psi(t))^T \Big[ -BR^{-1} B^T \Big]  
(P(t) x+ \Psi(t))  \\
+ (P(t) x+ \Psi(t))^T (AX +C) 
+ \frac{1}{2} \Big(X(t)^T  Q X(t) 
 +  u(t)^T R u(t)^T\Big) \\ 
  + L^T X(t) =0 
\   \mbox{ in } \mathcal S\times[0,T[,
\\ \\
P(T) =\phi, \quad \Psi(T) =0, \quad \chi(T)=0.
\end{array}\right.
\end{equation}

The boundary conditions are obtained by imposing that $$v(x,T)= \frac{1}{2} x^T P(T) x+ \Psi(T) x + \chi(T)=\frac{1}{2}x^T \phi x.$$

Since (\ref{cb33}) is an identity in $x$, it reduces to three equations:
\begin{equation}\label{cb4}
\left\{
\begin{array}{l}
\dot P + PA(x) + A(x)^T P - P BR^{-1} B^T P + Q=0
\ \mbox{in }  [0,T[,\ P(T) =\phi, \\ \\
\dot \Psi + A(x)^T \Psi + P C -  P BR^{-1} B^T \Psi + L =0 \ \mbox{in } [0,T[,\ \Psi(T) =0, \\ \\
\dot \chi  + \Psi^T C - \frac{1}{2} \Psi^T B R^{-1} B^T \Psi =0  \ \mbox{in } [0,T[,\ \chi(T) =0.
\end{array}\right.
\end{equation}

To understand the influence of the congestion term on the value function, let us now develop the expression for $\Psi$ and obtain
%

\begin{equation}
\begin{array}{lll}
\left[\begin{array}{c} \dot \Psi_1 \\ \dot \Psi_2\end{array}\right] 
+ \left[\begin{array}{ccc}
-\beta & 0 \\
k(x(t)) & 0 \end{array}\right]     
\left[\begin{array}{c}  \Psi_1 \\ 
 \Psi_2\end{array}\right] +
 \left[\begin{array}{cc}  P_{11} & P_{12} \\
P_{21} & P_{22}\end{array}\right] 
\left[\begin{array}{c}  \beta x_{off} \\ 0\end{array}\right] \\ \\ \qquad \qquad 
- \left[\begin{array}{c}  P_{12} (r_{on}^{-1} + r_{off}^{-1}) \Psi_2 \\
P_{22} (r_{on}^{-1} + r_{off}^{-1}) \Psi_2\end{array}\right] + 
\left[\begin{array}{c}  0 \\ Se+W\end{array}\right].\end{array}\end{equation}
The expression of $\Psi$ then can be rewritten as
\begin{equation}\label{cb10}
\displaystyle{ 
\left\{
\begin{array}{lll}
 \dot \Psi_1 - \beta \Psi_1 + P_{11} \beta x_{off} 
 - P_{12} (r_{on}^{-1} + r_{off}^{-1}) \Psi_2 =0,\\ \\
 \dot \Psi_2 + k(x(t)) \Psi_1 - P_{22} (r_{on}^{-1} + r_{off}^{-1}) \Psi_2 + (Se+W)=0,
\end{array}\right.
}
\end{equation}

which is of the form 
\begin{equation}\label{cb10}
\displaystyle{ 
\left\{
\begin{array}{lll}
 \dot \Psi_1 + a \Psi_1 + b \Psi_2 + c =0,\\ \\
\dot \Psi_2 + a' \Psi_1 + b' \Psi_2 + c' =0.
\end{array}\right.
}
\end{equation}
From the above set of inequalities, we obtain the solution $\Psi(x(t),e(t),t)$. Note that the term $a'$ depends on $x$ and $c'$ depends on $e(t)$.

Substituting the expression of the mean-field equilibrium strategies $u^* = - R^{-1} B^T [P X + \Psi]$ as in (\ref{mfes}) in the open-loop microscopic dynamics $\dot X(t) = A X(t) +  B u(t) + C$ introduced in (\ref{pq}), 
 and averaging both LHS and RHS we obtain the following closed-loop macroscopic dynamics
$$\dot{\bar X}(t) = [A(x) - B R^{-1} B^T P ] \bar X(t) - B R^{-1} B^T \bar \Psi(t)  + C,$$ 
where $\bar \Psi(t)=\int_{x_{on}}^{x_{off}} \int_{[0,1]} \Psi(x,e,t) m(x,y,t)dxdy$ and this concludes our proof. 

\subsection*{Proof of Corollary \ref{cor1}}

 Let $X(t)$ be a solution of dynamics (\ref{asdyn}) with initial value $X(0) \not \in \mathcal  X$. Set $t=\{\inf t >0|\, X(t) \in \mathcal X \}\leq \infty$. For all $t \in [0,t]$  
 \begin{equation}\nonumber
\begin{array}{lll}
V(X(t+ d t)) - V(X(t)) & = & \| X(t+ d t) - \Pi_{\mathcal X}(X(t))\| - \|X(t) -\Pi_{\mathcal X}(X(t))\| \\ 
& = & \| X(t) + d X(t) - \Pi_{\mathcal X}(X(t))\| - \| X(t) - \Pi_{\mathcal X}(X(t))\| 
\\ 
& = &\frac{1}{\| X(t) + d X(t) - \Pi_{\mathcal X}(X(t))\|} \| X(t) + d X(t) - \Pi_{\mathcal X}(X(t))\|^2 \\
&&- \frac{1}{\| X(t) - \Pi_{\mathcal X}(X(t))\|} \| X(t) - \Pi_{\mathcal X}(X(t))\|^2. 
\end{array}\end{equation}
Taking the limit of the difference above we obtain
\begin{equation}\nonumber
\displaystyle \begin{array}{lll}
\dot V (X(t)) & = & \lim_{dt \rightarrow 0} \frac{ V(X(t+dt)) - V (X(t))}{dt}\\ 
& = &\lim_{dt \rightarrow 0}  \frac{1}{dt} \Big[\frac{1}{\| X(t) + d X(t) - \Pi_{\mathcal X}(X(t))\|} \| X(t) + d X(t) - \Pi_{\mathcal X}(X(t))\|^2\\ && - \frac{1}{\| X(t) - \Pi_{\mathcal X}(X(t))\|} \|X(t) - \Pi_{\mathcal X}(X(t))\|^2 \Big]\\
& \leq & \frac{1}{\|X(t) - \Pi_{\mathcal X}(X(t))\|} \Big[  \partial_X V(X,t)^T \Big([A - B R^{-1} B^T  P] X(t) \\ && - B R^{-1} B^T \Psi^*(x(t),e(t))  + C\Big) 
 + \|X(t) - \Pi_{\mathcal X}(X(t))\|^2\Big] < 0, 
\end{array}\end{equation}
which implies $\mathcal L V (X(t)) <0$, for all $X(t) \not \in \mathcal X$ and this concludes our proof. 

\subsection*{Proof of Theorem \ref{thm2}}
This proof follows the same reasoning as the proof of Theorem \ref{thm1}. 
However, differently from there, here for the quadratic terms in (\ref{cb333}) we have
$$ \sigma_{11}(x)^2 P_{11}(t)+\sigma_{22}(y)^2 P_{22}(t)= \hat \sigma_{11}^2 x^2 P_{11}(t)+\hat \sigma_{22}^2 y^2 P_{22}(t).$$
Reviewing (\ref{cb333})  as an identity in $x$, this leads  to the following three equations to solve in the variable $P(t)$, $\Psi(t)$, and $\chi(t)$:
\begin{equation}\label{cb4}
\left\{
\begin{array}{l}
\dot P(t)  + P(t)A + A^T P - P BR^{-1} B^T P + Q+ \tilde P=0
\ \mbox{in }  [0,T[,\ P(T) =\phi, \\ \\
\dot \Psi(t) + A^T \Psi + P C -  P BR^{-1} B^T \Psi + L =0 \ \mbox{in } [0,T[,\ \Psi(T) =0, \\ \\
\dot \chi(t)  + \Psi(t)^T C - \frac{1}{2} \Psi^T B R^{-1} B^T \Psi =0  \ \mbox{in } [0,T[,\ \chi(T) =0,
\end{array}\right.
\end{equation}
where 
\begin{equation}\label{diag}
\tilde P = Diag( (\hat \sigma_{ii}^2  P_{ii})_{i =1,2} )=
\left[\begin{array}{cc} \hat \sigma_{11}^2  P_{11} & 0 \\
0 & \hat \sigma_{22}^2 P_{22} \end{array}\right].\end{equation} 

\subsection*{Proof of Corollary \ref{cor2}}
 Let $X(t)$ be a solution of dynamics (\ref{asdyn1111}) with initial value $X(0) \not \in \mathcal X$. Set $t=\{\inf t >0|\, X(t) \in \mathcal X \}\leq \infty$ and let  $V(X(t)) = dist(X(t),\mathcal X)$.   For all $t \in [0,t]$  
\begin{equation}\nonumber
\begin{array}{lll}
V(X(t+ d t)) - V(X(t)) & = & \| X(t+ d t) - \Pi_{\mathcal X}(X(t))\| - \|X(t) -\Pi_{\mathcal X}(X(t))\| \\ 
& = & \| X(t) + d X(t) - \Pi_{\mathcal X}(X(t))\| - \| X(t) - \Pi_{\mathcal X}(X(t))\| 
\\ 
& = &\frac{1}{\| X(t) + d X(t) - \Pi_{\mathcal X}(X(t))\|} \| X(t) + d X(t) - \Pi_{\mathcal X}(X(t))\|^2 - \\
&& \frac{1}{\| X(t)  - \Pi_{\mathcal X}(X(t))\|} \| X(t) - \Pi_{\mathcal X}(X(t))\|^2. 
\end{array}\end{equation}
From the definition of infinitesimal generator 
\begin{equation}\nonumber
\displaystyle \begin{array}{lll}
\mathcal L V (X(t)) & = & \lim_{dt \rightarrow 0} \frac{\mathbb E V(X(t+dt)) - V (X(t))}{dt}\\
& = &\lim_{dt \rightarrow 0}  \frac{1}{dt} \Big[\mathbb E \Big(\frac{1}{\| X(t) + d X(t) - \Pi_{\mathcal X}(X(t))\|} \| X(t) + d X(t) -\Pi_{\mathcal X}(X(t))\|^2\Big)\\ && - \frac{1}{\| X(t) - \Pi_{\mathcal X}(X(t))\|} \|X(t) - \Pi_{\mathcal X}(X(t))\|^2 \Big]\\ 
& \leq & \frac{1}{\|X(t) - \Pi_{\mathcal X}(X(t))\|} 
\Big[  \partial_X V(X,t)^T \Big([A - B R^{-1} B^T  P] X(t) \\ && - B R^{-1} B^T \Psi^*(x(t),e(t))  + C\Big) 
  \\  && \qquad \qquad  \qquad  + \frac{1}{2} ( \sigma_{11}^2(x)  \partial_{xx} V(X,t)
 +  \sigma_{22}^2(y)  \partial_{yy} V(X,t) )\Big].
\end{array}\end{equation}

From (\ref{contract}) the above implies that 
 $\mathcal L V (X(t)) <0$, for all $X(t) \not \in  \mathcal M$ and this concludes our proof. 

\subsection*{Proof of Theorem \ref{thm3}}
From (\ref{case2}),  in the HJB equation (\ref{cb333})  we now have constant terms
$$\frac{1}{2} \sum_{i=1}^2  \sigma_{ii}(.)^2 P_{ii}(t)= \hat \sigma_{11}^2  P_{11}(t)+\hat \sigma_{22}^2  P_{22}(t).$$
Again, since  the HJB equation (\ref{cb333})  is an identity in $x$, it reduces to three equations:
\begin{equation}\label{cb4}
\left\{
\begin{array}{l}
\dot P(t)  + P(t)A + A^T P - P BR^{-1} B^T P + Q=0
\ \mbox{in }  [0,T[,\ P(T) =\phi, \\ \\
\dot \Psi(t) + A^T \Psi + P C -  P BR^{-1} B^T \Psi + L =0 \ \mbox{in } [0,T[,\ \Psi(T) =0, \\ \\
\dot \chi(t)  + \Psi(t)^T C - \frac{1}{2} \Psi^T B R^{-1} B^T \Psi + \tilde P =0  \ \mbox{in } [0,T[,\ \chi(T) =0,
\end{array}\right.
\end{equation}
where 
\begin{equation}\label{diag}
\tilde P = 
\left[\begin{array}{cc} \hat \sigma_{11}^2& 0 \\
0 & \hat \sigma_{22}^2  \end{array}\right].\end{equation}

Substituting the expression of the mean-field equilibrium strategy $u^* = - R^{-1} B^T [P X + \Psi]$ as in (\ref{mfes1}) in the open-loop microscopic dynamics $d X(t) = (A X(t) +  B u(t) + C) dt + \Sigma d\mathcal B_t$ given in (\ref{pq1})
and averaging both LHS and RHS we obtain the following closed-loop macroscopic dynamics
$$\dot{\bar X}(t) = [A - B R^{-1} B^T P ] \bar X(t) - B R^{-1} B^T \Psi^*(\bar X(t))  + C,$$
and this concludes our proof.

\subsection{Proof of Corollary \ref{cor3}}
 Let $X(t)$ be a solution of dynamics (\ref{asdyn11}) with initial value $X(0) \not \in \mathcal X$. Set $t=\{\inf t >0|\, X(t) \in \mathcal X \}\leq \infty$ and let  $V(X(t)) = dist(X(t),\mathcal X)$.   For all $t \in [0,t]$  
\begin{equation}\nonumber
\begin{array}{lll}
V(X(t+ d t)) - V(X(t)) & = & \| X(t+ d t) - \Pi_{\mathcal X}(X(t))\| - \|X(t) -\Pi_{\mathcal X}(X(t))\| \\ 
& = & \| X(t) + d X(t) - \Pi_{\mathcal X}(X(t))\| - \| X(t) - \Pi_{\mathcal X}(X(t))\| 
\\ 
& = &\frac{1}{\| X(t) + d X(t) - \Pi_{\mathcal X}(X(t))\|} \| X(t) + d X(t) - \Pi_{\mathcal X}(X(t))\|^2 - \\
&& \frac{1}{\| X(t)  - \Pi_{\mathcal X}(X(t))\|} \| X(t) - \Pi_{\mathcal X}(X(t))\|^2 
\end{array}\end{equation}
From the definition of infinitesimal generator 
\begin{equation}\nonumber
\displaystyle \begin{array}{lll}
\mathcal L V (X(t)) & = & \lim_{dt \rightarrow 0} \frac{\mathbb E V(X(t+dt)) - V (X(t))}{dt}\\
& = &\lim_{dt \rightarrow 0}  \frac{1}{dt} \Big[\mathbb E \Big(\frac{1}{\| X(t) + d X(t) - \Pi_{\mathcal X}(X(t))\|} \| X(t) + d X(t) -\Pi_{\mathcal X}(X(t))\|^2\Big)\\ && - \frac{1}{\| X(t) - \Pi_{\mathcal X}(X(t))\|} \|X(t) - \Pi_{\mathcal X}(X(t))\|^2 \Big]\\ 
& \leq & \frac{1}{\|X(t) - \Pi_{\mathcal X}(X(t))\|} 
\Big[  \partial_X V(X,t)^T \Big([A - B R^{-1} B^T  P] X(t) \\ && - B R^{-1} B^T \Psi^*(x(t),e(t))  + C\Big) 
  \\  && \qquad \qquad  \qquad  + \frac{1}{2} ( \hat \sigma_{11}^2  \partial_{xx} V(X,t)
 + \hat \sigma_{22}^2  \partial_{yy} V(X,t) )\Big].
\end{array}\end{equation}

From (\ref{contract1}) the above implies that 
 $\mathcal L V (X(t)) <0$, for all $X(t) \not \in  \mathcal M$ and this concludes our proof. 

\subsection*{Proof of Theorem \ref{thm4}}

Isolating the HJI equation in (\ref{mainHJI}), we have
\begin{equation}\label{1111expla}\nonumber
\left\{
\begin{array}{lll}
-\partial_t \mathcal V_t(X) - \sup_u \inf_w\Big\{  \partial_X \mathcal V_t(X)^T (AX + Bu +C+D w) + \frac{1}{2} \Big(X(t)^T  Q X(t) \\
 \qquad \qquad +  u(t)^T R u(t) - \gamma^2 w(t)^T w(t)  \Big)+ L^T X(t) \Big\}=0, \mbox{ in } \mathcal S \times [0,T[,\\ \\
\mathcal V_T(X)=g(x) \mbox{ in } \mathcal S.  
\end{array}\right.
\end{equation}

\noindent

Let us consider the following value function 
$$v(X,t)= \frac{1}{2}X^TP(t) X+ \Psi(t)^T X + \chi(t),$$ 

and the corresponding mean-field equilibrium control and worst-case disturbance 
\begin{equation}\nonumber
\begin{array}{ccl}
u^* &=& - R^{-1} B^T [P X + \Psi],\\
w^* &=& \frac{1}{\gamma^2} D^T [P X + \Psi].
\end{array}
\end{equation}

so that (\ref{1111expla}) can be rewritten as

\begin{equation}\label{cb3}\displaystyle
\left\{
\begin{array}{r}
\frac{1}{2} X^T\dot P(t) X + \dot \Psi(t)X + \dot \chi(t) +(P(t) X+ \Psi(t))^T \Big[ -BR^{-1} B^T
+ \frac{1}{\gamma^2} D D^T \Big]  
(P(t) x+ \Psi(t))  \\
+ (P(t) x+ \Psi(t))^T (AX +C) 
+ \frac{1}{2} \Big(X(t)^T  Q X(t) 
 +  u(t)^T R u(t) - \gamma^2 w(t)^T w(t)\Big) \\ 
  + L^T X(t) +\frac{1}{2} \sum_{i=1}^2  \sigma_{ii}(.)^2 P_{ii}(t)=0 
\ \mbox{in } \mathbb{R}^2\times[0,T[,
\\ \\
P(T) =\phi, \quad \Psi(T) =0, \quad \chi(T)=0.
\end{array}\right.
\end{equation}

The boundary conditions are obtained by imposing that $$v(X,T)= \frac{1}{2} X^T P(T) X+ \Psi(T) X + \chi(T)=\frac{1}{2}X^T \phi X.$$


%
The above set of identities  in $x$ yields the following three equations in the variable $P(t)$, $\Psi(t)$, and $\chi(t)$:
\begin{equation}\label{cbwc}
\left\{
\begin{array}{l}
\dot P(t)  + P(t)A + A^T P +P (- B R^{-1} B^T  + \frac{1}{\gamma^2} DD^T) P + Q=0
\ \mbox{in }  [0,T[,\ P(T) =\phi, \\ \\
\dot \Psi(t) + A^T \Psi + P C +(- B R^{-1} B^T  + \frac{1}{\gamma^2} DD^T) \Psi + L =0 \ \mbox{in } [0,T[,\ \Psi(T) =0, \\ \\
\dot \chi(t)  + \Psi(t)^T C + \frac{1}{2} \Psi^T (- B R^{-1} B^T  + \frac{1}{\gamma^2} DD^T) \Psi =0  \ \mbox{in } [0,T[,\ \chi(T) =0.
\end{array}\right.
\end{equation}

Substituting the expressions of the mean-field equilibrium strategies $u^* = - R^{-1} B^T [P X + \Psi]$ and $w^* = \frac{1}{\gamma^2} D^T [P X + \Psi]$ as in (\ref{mfes4}) in the open-loop microscopic dynamics $\dot X(t) = A X(t) +  B u(t) + C$ introduced in (\ref{pq4}), 
 and averaging both LHS and RHS we obtain the following closed-loop macroscopic dynamics
$$\dot{\bar X}(t) = [A +(- B R^{-1} B^T  + \frac{1}{\gamma^2} DD^T)P] \bar X(t) +(- B R^{-1} B^T +  \frac{1}{\gamma^2} DD^T) \Psi^*(\bar X(t))  + C,$$ and this concludes our proof.  

\subsection*{Proof of Corollary \ref{cor4}}

 Let $X(t)$ be a solution of dynamics (\ref{wcdyn}) with initial value $X(0) \not = \mathcal X$. Set $t=\{\inf t >0|\, X(t) \in \mathcal X \}\leq \infty$ and let  $V(X(t)) = dist(X(t),\mathcal X)$.   For all $t \in [0,t]$  
\begin{equation}\nonumber
\begin{array}{lll}
V(X(t+ d t)) - V(X(t)) & = & \| X(t+ d t) - \Pi_{\mathcal X}(X(t))\| - \|X(t) -\Pi_{\mathcal X}(X(t))\| \\ 
& = & \| X(t) + d X(t) - \Pi_{\mathcal X}(X(t))\| - \| X(t) - \Pi_{\mathcal X}(X(t))\| 
\\ 
& = &\frac{1}{\| X(t) + d X(t) - \Pi_{\mathcal X}(X(t))\|} \| X(t) + d X(t) - \Pi_{\mathcal X}(X(t))\|^2\\
&&  - \frac{1}{\| X(t) - \Pi_{\mathcal X}(X(t))\|} \| X(t) - \Pi_{\mathcal X}(X(t))\|^2. 
\end{array}\end{equation}
From the definition of infinitesimal generator 
\begin{equation}\nonumber
\displaystyle \begin{array}{lll}
\dot V (X(t)) & = & \lim_{dt \rightarrow 0} \frac{ V(X(t+dt)) - V (X(t))}{dt}\\ 
& = &\lim_{dt \rightarrow 0}  \frac{1}{dt} \Big[\frac{1}{\| X(t) + d X(t) - \Pi_{\mathcal X}(X(t))\|} \| X(t) + d X(t) - \Pi_{\mathcal X}(X(t))\|^2\\ 
& \leq & \frac{1}{\|X(t) -  \Pi_{\mathcal X}(X(t))\|} \Big[ \partial_X V(X,t)^T \Big([A +(- B R^{-1} B^T  + \frac{1}{\gamma^2} DD^T)P] X(t) \\ && +(- B R^{-1} B^T +  \frac{1}{\gamma^2} DD^T)\Psi^*(x(t),e(t))  + C\Big)
\leq 0 
\end{array}\end{equation}

%
%

which implies $\mathcal L V (\rho(t)) <0$, for all $X(t) \not = \mathcal X$ and this concludes our proof. 

\end{document}